\documentclass[12pt,reqno]{amsart}                  
\usepackage{amssymb}
\usepackage{mathrsfs}

\newtheorem{theorem}{Theorem}
\newtheorem{lemma}{Lemma}
\begin{document}
\title[On the dynamic programming principle]{On the dynamic programming principle for controlled diffusion processes in a cylindrical region}
\author{Dmitry B. Rokhlin}

\address{D.B. Rokhlin,
Faculty of Mathematics, Mechanics and Computer Sciences,
              Southern Federal University, 
Mil'chakova str., 8a, 344090, Rostov-on-Don, Russia}          
\email{rokhlin@math.rsu.ru}                   

\begin{abstract} We prove the dynamic programming principle for a class of diffusion processes controlled up to the time of exit from a cylindrical region $[0,T)\times G$. It is assumed that the functional to be maximized is in the Lagrange form with nonnegative integrand. Besides this we only adopt the standard assumptions, ensuring the existence of a unique strong solution of a stochastic differential equation for the state process.
\end{abstract}

\maketitle

\section{Main result} 
\setcounter{equation}{0}
Bellman's dynamic programming (or optimality) principle is a fundamental result of optimal control. In diffusion setup it connects a stochastic optimal control problem with a Hamilton-Jacobi-Bellman partial differential equation \cite{Kry77}, \cite{FleSon06}, \cite{Tou13}. 

In this note we prove the validity of this principle for a diffusion process controlled up to the time of exit from a cylindrical region $[0,T)\times G$, where $G\subset\mathbb R^d$ is an open set. We consider the maximization problem for a functional in the Lagrange form with nonnegative integrand. This case was also isolated in \cite{DonKry07}, where the justification of the dynamic programming principle was based on the reduction to the problem with optimal stopping and then to the case $G=\mathbb R^d$ (as it was done in \cite{Kry81}). The proof presented below is a more direct one. In fact, it is shown that the argumentation of \cite{BouTou11} can be adapted to the present case. It should be mentioned that we use several ideas of the papers \cite{BayHua10}, \cite{EsfChaLyg12}. Note also that the case of infinite time horizon is not addressed in the present paper. 

Throughout the paper we work on the canonical probability space. Namely, let $\Omega=C:=C([0,T];\mathbb R^d)$ be the Banach space of continuous $\mathbb R^d$-valued functions with the norm $\|\omega\|=\sup_{0\le t\le T} |\omega(t)|$ and let $\mathsf P$ be the Wiener measure on its Borel $\sigma$-algebra. So, the canonical process $W_s(\omega)=\omega_s$ is the standard $d$-dimensional Brownian motion under $\mathsf P$. Denote by $\mathbb F=(\mathscr F_t)_{0\le t\le T}$
the minimal augmented filtration (see e.g. \cite{Bas11}, Chapter 20) generated by the coordinate mappings $W_s(\omega)=\omega_s$, $0\le s\le T$.

Furthermore, for each $t\in[0,T]$ we introduce the filtration $\mathbb F^t=(\mathscr F^t_s)_{0\le s\le T}$, where $\mathscr F^t_s$ is generated by the increments $W_v-W_u$, $t\le u\le v\le s\vee t$ and $\mathscr F_0$.
The filtration $\mathbb F^t$ represents the information flow if the observation process starts at time $t$. The information contained in $\mathscr F_0$ (and concerning $\mathsf P$-null sets) is available at any time. We have $\mathbb F^0=\mathbb F$. 

Let $U$ be a separable metric space and assume that $U=\bigcup_{n=1}^\infty U(n)$ for some sequence $U(n)\subseteq U(n+1)$ of its subsets. Denote by $\mathscr A^t(n)$ the set of $\mathbb F^t$-progressively measurable processes with values in $U(n)$ and put $\mathscr A^t=\bigcup_{n=1}^\infty \mathscr A^t(n)$. Given a point $(t,x)\in [0,T]\times\mathbb R^n$ and a control process $\alpha\in\mathscr A:=\mathscr A^0$, consider the stochastic differential equation
\begin{equation} \label{eq:1.1}
dX_s^{t,x,\alpha}=b(s,X_s^{t,x,\alpha},\alpha_s)ds+\sigma(s,X_s^{t,x,\alpha},\alpha_s)dW_s,\ \ \ s\in [t,T]
\end{equation}
for the state process $X^{t,x,\alpha}$ with the initial condition $X_t^{t,x,\alpha}=x$. 

We assume that the functions
$$ b:[0,T]\times\mathbb R^d\times U\mapsto\mathbb R^d, \ \ \
   \sigma:[0,T]\times\mathbb R^d\times U\mapsto\mathbb R^d\times\mathbb R^d$$
are Borel and satisfy the conditions
\begin{subequations}
\begin{align*}
|b(t,x,u)-b(t,y,u)|+|\sigma(t,x,u)-\sigma(t,y,u)| & \le K_n |x-y| \\
|b(t,x,u)|+|\sigma(t,x,u)| & \le K_n (1+|x|) 
\end{align*}
\end{subequations}
for $(t,x)\in [0,T]\times\mathbb R^d$ and $u,v\in U(n)$ with some constants $K_n$, possibly depending on $n$. 
Under these assumptions there exists a unique strong solution  $X^{t,x,\alpha}$ of (\ref{eq:1.1}) for any $\alpha\in\mathscr A$. We put  $X^{t,x,\alpha}_s=x$ for $s\in [0,t]$. 

Fix an open set $G\subset\mathbb R^d$ and put
$$\tau^{t,x,\alpha}=\inf\{s\ge t:X_s^{t,x,\alpha}\not\in G\}\wedge T.$$ We consider the stochastic optimal control problem 
$$ v(t,x):=\sup_{\alpha\in\mathscr A^t} J(t,x,\alpha):=\sup_{\alpha\in\mathscr A^t}\mathsf E\int_{t}^{\tau^{t,x,\alpha}} f(s,X^{t,x,\alpha}_s,\alpha_s)\,ds, $$
where $f$ is a \emph{nonnegative} Borel function on $[0,T]\times\mathbb R^d\times U$ and $f(t,\cdot,u)$ is lower semicontinuous on $\mathbb R^d$. 

\begin{theorem} Let $\mathcal{T}^t_{t,\tau}$ be the set of all $\mathbb F^t$-stopping times $\theta$ such that $t\le\theta\le\tau$. Then the function $v$ is lower semicontinuous on $[0,T]\times\mathbb R^d$ and
\begin{align*}
v(t,x) &=\sup_{\alpha\in\mathscr A^t}\inf_{\theta\in \mathcal T^t_{t,\tau^{t,x,\alpha}}}
\mathsf E\left[\int_t^\theta f(s,X_s^{t,x,\alpha},\alpha_s)\,ds+v(\theta,X_\theta^{t,x,\alpha})\right],\\
 &=\sup_{\alpha\in\mathscr A^t}\sup_{\theta\in \mathcal T^t_{t,\tau^{t,x,\alpha}}}
\mathsf E\left[\int_t^\theta f(s,X_s^{t,x,\alpha},\alpha_s)\,ds+v(\theta,X_\theta^{t,x,\alpha})\right].
\end{align*}
\end{theorem}

A similar result in a more general setting was proved in \cite{DonKry07} (Theorem 2.1). However, we do not assume that the functions $b$, $\sigma$, $f$ are bounded and continuous in control variable and that the domain $G$ is bounded. For the objective functional in the Mayer form an assertion close to Theorem 1 was proved in \cite{EsfChaLyg12} (Theorem 4.7)  under the assumption that the diffusion term is non-degenerate. We mention also the result of \cite{FleSon06} (Chapter V, Theorem 2.1), where the dynamic programming principle was proved under the assumptions, ensuring the continuity of the value function $v$. Other references can be found in the cited literature.

In Section 2 we collect some auxiliary results. None of them is new. Theorem 1 is proved in Section 3 with the use of
the technique of \cite{BouTou11} adapted to the exit time problem. Similar approach was used in \cite{EsfChaLyg12}. 

\section{Auxilary results}
\setcounter{equation}{0}
The following lemma states the existence of  a continuous modification $Y^{t,x,\alpha}$ of the multiparameter process $X^{t,x,\alpha}$ (see \cite{Bor89}, Theorem 2.1). This fact is a consequence  of Kolmogorov's continuity theorem \cite{RogWil00} (Theorem 25.2) or \cite{Pro04} (Chapter 4, Theorem 72). We mention that it was used in \cite{EsfChaLyg12} (in the course of the proof of Proposition 4.3). For reader's convenience, we outline the proof of this result. Some details are borrowed from the lectures on stochastic analysis by Timo Sepp\"al\"ainen.   
\begin{lemma} \label{lem:1} Fix $\alpha\in\mathscr A$. There exists an $\mathbb F$-adapted process such that
\begin{itemize}
\item[(i)] $Y^{t,x,\alpha}_\cdot(\omega)\in C$ and the mappings 
$$(t,x)\mapsto Y^{t,x,\alpha}_\cdot(\omega): [0,T]\times\mathbb R^d\mapsto C$$
are continuous for each $\omega\in\Omega$;
\item[(ii)] $\|Y^{t,x,\alpha}_\cdot(\omega)-X^{t,x,\alpha}_\cdot(\omega)\|_C=0$ $\mathsf P$-a.s. for all $(t,x)\in [0,T]\times\mathbb R^d$.
\end{itemize}
\end{lemma}
\emph{Proof.} 
Denote by $X^{t,\xi,\alpha}$ the solution of (\ref{eq:1.1}) with the initial condition $X_t^{t,\xi,\alpha}=\xi\in L^{2q}(\mathscr F_t)$ (note that for random $\xi$ we do not define the process $X^{t,\xi,\alpha}_u$ for $u<t$). For brevity let's omit index $\alpha$ till the end of the proof. As is known (see \cite{Kry77}, Corollary II.5.10 and Theorem II.5.9), for $q\ge 1$ the following inequalities hold true:
\begin{align} 
\mathsf E\sup_{t\le u\le s}|X_u^{t,\xi}-\xi|^{2q} \le & K (s-t)^q(1+\mathsf E|\xi|^{2q}), \label{eq:2.1}\\
\mathsf E\sup_{t\le u\le T}|X_u^{t,\xi}-X_u^{t,\eta}|^{2q} \le & K\mathsf E|\xi-\eta|^{2q}. \label{eq:2.2}
\end{align}
Let $s>t$. From the inequality 
\begin{align*}
|X_u^{s,y}-X_u^{t,x}|\le |X_u^{s,y}-X_u^{t,y}|+|X_u^{t,y}-X_u^{t,x}|=|y-X_u^{t,y}| I_{[t,s]}(u)\\
+|X_u^{s,y}-X_u^{s,X_s^{t,y}}| I_{[s,T]}(u)+|x-y|I_{[0,t]}(u)+ |X_u^{t,y}-X_u^{t,x}| I_{[t,T]}(u),
\end{align*}
where the flow property $X_u^{t,y}=X_u^{s,X_s^{t,y}}$, $u\ge s$ ($\mathsf P$-a.s.) is used (see e.g. \cite{BayHua10}, Remark 2.5(i)), it follows that
\begin{align*}
\sup_{0\le u\le T} |X_u^{s,y}-X_u^{t,x}|^{2q}\le  4^{2q-1} \left(\sup_{t\le u\le s} |y-X_u^{t,y}|^{2q}+\right.\\
\left. +\sup_{s\le u\le T}  |X_u^{s,y}-X_u^{s,X_s^{t,y}}|^{2q} +|x-y|^{2q}+  \sup_{t\le u\le T} |X_u^{t,y}-X_u^{t,x}|^{2q}\right).
\end{align*}
Estimating according to (\ref{eq:2.1}), (\ref{eq:2.2}) the expectations of the terms on the right-hand side of the last inequality:
\begin{align*}
\mathsf E \sup_{t\le u\le s} |y-X_u^{t,y}|^{2q}\le K (s-t)^q(1+|y|^{2q}) \\
\mathsf E \sup_{s\le u\le t} |X_u^{s,y}-X_u^{s,X_s^{t,y}}|^{2q}\le K \mathsf E |y-X_s^{t,y}|^{2q}\le K^2 (s-t)^q(1+|y|^{2q})\\
\mathsf E \sup_{t\le u\le T} |X_u^{t,y}-X_u^{t,x}|^{2q}\le K |y-x|^{2q}, 
\end{align*}
and interchanging $s$ and $t$, we get
\begin{align} 
\mathsf E \sup_{0\le u\le T} |X_u^{s,y}-X_u^{t,x}|^{2q} &\le K'( |s-t|^q (1+|y|^{2q})+|x-y|^{2q}) \nonumber\\
&\le K_n' (|s-t|^2+|y-x|^2)^{q/2} \label{eq:2.3}
\end{align}
for $(t,x)\in [0,T]\times (n \Delta)$, where $\Delta$ is the hypercube $[-1,1]^d$ in $\mathbb R^d$. 

Put $D_k=2^{-k}\mathbb Z^{d+1}$ and denote by $D=\cup_{k=1}^\infty D_k$ the set of dyadic rational points in $\mathbb R^{d+1}$. The estimate (\ref{eq:2.3}) implies that the restriction of $X^{t,x}_\cdot$ to $([0,T]\times (n \Delta)) \cap D$ is a.s. uniformly continuous as a function with values in $C$ (see \cite{Pro04}, Chap. 4, the proof of Theorem 72). Thus, there exists a set $\Omega'\subset\mathscr F$, $\mathsf P(\Omega')=1$ such that the function 
$$(t,x)\mapsto X^{t,x}_\cdot(\omega): ([0,T]\times \mathbb R^d) \cap D\mapsto C,\ \ \omega\in\Omega'$$
is uniformly continuous on compact sets. This function is uniquely extendable to the continuous function $Y^{t,x}_\cdot(\omega)$ on $[0,T]\times \mathbb R^d$. For $\omega\not\in\Omega'$ we put $Y^{t,x}_s(\omega)=x$, $s\in[0,T]$.

It remains to check (ii) for $(t,x)\not\in D$. Take a sequence $(t_n,x_n)\to (t,x)$, $(t_n,x_n)\in ([0,T]\times\mathbb R^d)\cap D$.
We have
\begin{align*}
\mathsf P(\|X^{t,x}-Y^{t,x}\|_C\ge\varepsilon)\le \mathsf P(\|X^{t,x}-X^{t_n,x_n}\|_C\ge\varepsilon/3)\\
+ \mathsf P(\|X^{t_n,x_n}-Y^{t_n,x_n}\|_C\ge\varepsilon/3)+\mathsf P(\|Y^{t_n,x_n}-Y^{t,x}\|_C\ge\varepsilon/3)
\end{align*}
As $n\to \infty$, the first term on the right vanishes by the Chebyshev inequality and (\ref{eq:2.3}), the third term vanishes by the continuity of $Y^{t,x}$ and the second term is equal to zero by the definition of $Y^{t,x}$. It follows that 
$\mathsf P(\|X^{t,x}-Y^{t,x}\|_C>0)=0$. $\square$  

For $X\in C$ and an open set $G\subset\mathbb R^d$ put
$$ \tau_G(X)=\inf\{0\le t\le T: X(t)\not\in G\}\wedge T.$$ The following result is of folklore type. See however \cite{Aub09} (Lemma 4.2.2) for the case of infinite time horizon.
\begin{lemma} \label{lem:2} The function $\tau_G: C\mapsto \mathbb R$ is lower semicontinuous in the norm topology of $C$.
\end{lemma}

\emph{Proof}. We show that the set $L_a=\{X\in C:\tau_G(X)>a\}$ is open in $C$ for all $a\in\mathbb R$. Note that $L_a=C$ for $a<0$. Let $X\in L_a$ with $a\ge 0$. Then $X(t)\in G$ for all $t$ in the compact interval $[0,a]\subset[0,T]$. Denote by $G^c$ the compliment of $G$ in $\mathbb R^d$. Since the distance function
$$ \rho(x,A)=\inf\{|x-y|:y\in A\},\ \ \ A\subset\mathbb R^d$$
is continuous in $x\in\mathbb R^d$:
$$ |\rho(x,A)-\rho(y,A)|\le |x-y|$$
(see \cite{AliBor06}, Theorem 3.16) and $\rho(X(t),G^c)>0$, $t\in [0,a]$, we have
$$ \inf_{t\in [0,a]} \rho(X(t),G^c)=\delta>0.$$

Consider the neighborhood of $X$ 
$$U=\{Y\in C: \sup_{t\in [0,T]} |Y(t)-X(t)|<\delta/2\}.$$
The inequality
$$ \rho(Y(t),G^c)\ge\rho(X(t),G^c)-|X(t)-Y(t)|\ge \delta-\delta/2=\delta/2,\ \ t\in [0,a]$$
shows that $\tau_G(Y)>a$ for all $Y\in U$. $\square$

\begin{lemma} \label{lem:3} 
Fix $\alpha\in\mathscr A$. The function $(t,x)\mapsto\tau^{t,x,\alpha}(\omega)$ is lower semicontinuous on $[0,T]\times\mathbb R^d$ for $\mathsf P$-a.e. $\omega\in\Omega$.
\end{lemma}

\emph{Proof}. Note that for fixed $(t,x)$ the process $(X^{t,x,\alpha})_{0\le s\le T}$ and its continuous modification $(Y^{t,x,\alpha})_{0\le s\le T}$, constructed in Lemma \ref{lem:1}, are indistinguishable. Hence, 
$\tau^{t,x,\alpha}(\omega)$ is $\mathsf P$-a.s. a superposition of the continuous function
$$ (t,x)\mapsto Y_\cdot^{t,x,\alpha}(\omega):[0,T]\times\mathbb R^d\mapsto C$$
and the lower semicontinuous function $\tau_G:C\mapsto [0,T]$. $\qed$

\begin{lemma} \label{lem:4}
For any $\alpha\in\mathscr A$ the function 
$$ (t,x)\mapsto J(t,x,\alpha):[0,T]\times\mathbb R^d\mapsto [0,+\infty]$$
is lower semicontinuous.
\end{lemma}
\emph{Proof}. Put $B_\delta(t,x)=\{(s,y)\in [0,T]\times\mathbb R^d:|s-t|^2+|y-x|^2\le\delta^2\}$ and denote by
$$J_*(t,x,\alpha)=\lim_{n\to\infty}\inf_{(s,y)\in B_{1/n}(t,x)} J(s,y,\alpha),$$
the lower semicontinuous envelope of $J$. For $(t,x)\in [0,T]\times\mathbb R^d$ take a sequence $(t_n,x_n)\in B_{1/n}(t,x)$ such that
$$ J_*(t,x,\alpha)=\lim_{n\to\infty} J(t_n,x_n,\alpha).$$
By the Fatou lemma we have
$$ J_*(t,x,\alpha)\ge \mathsf E\liminf_{n\to\infty}\int_{t_n}^{\tau^{t_n,x_n,\alpha}} f(s,X^{t_n,x_n,\alpha}_s,\alpha_s)\,ds.$$
Furthermore, by the ($\mathsf P$-a.s.) lower semicontinuity of the function $(t,x)\mapsto\tau^{t,x,\alpha}(\omega)$ for any $\varepsilon>0$ the exists $N(\omega)$ such that  
$$ \tau^{t_n,x_n,\alpha}(\omega)\ge \tau^{t,x,\alpha}(\omega)-\varepsilon,\ \ n\ge N(\omega).$$
Since $f$ is nonnegative it follows that
$$ J_*(t,x,\alpha)\ge \mathsf E\liminf_{n\to\infty}\int_0^T I_{(t_n,\tau^{t,x,\alpha}-\varepsilon)}(s) f(s,X^{t_n,x_n,\alpha}_s,\alpha_s)\,ds.$$
By the continuity of $(t,x)\mapsto X^{t,x,\alpha}_s(\omega)$ and the lower semicontinuity of $x\mapsto f(t,x,a)$
we have
\begin{align*}
\liminf_{n\to\infty}I_{(t_n,\tau^{t,x,\alpha}-\varepsilon)}(s) f(s,X^{t_n,x_n,\alpha}_s,\alpha_s)\ge &
I_{(t,\tau^{t,x,\alpha}-\varepsilon)}\liminf_{n\to\infty}f(s,X^{t_n,x_n,\alpha}_s,\alpha_s)\\
\ge & I_{(t,\tau^{t,x,\alpha}-\varepsilon)}f(s,X^{t,x,\alpha}_s,\alpha_s).
\end{align*}
Thus, by the Fatou lemma,
$$ J_*(t,x,\alpha)\ge \mathsf E \int_0^T I_{(t,\tau^{t,x,\alpha}-\varepsilon)}f(s,X^{t,x,\alpha}_s,\alpha_s)\,ds.$$
At last, by the monotone convergence theorem we get
$$ J_*(t,x,\alpha)\ge\lim_{\varepsilon\downarrow 0 }\mathsf E \int_t^{\tau^{t,x,\alpha}-\varepsilon}f(s,X^{t,x,\alpha}_s,\alpha_s)\,ds=J(t,x,\alpha).$$
The converse inequality is evident. $\square$ 

For an $\mathbb F$-stopping time $\theta\in\mathcal T:=\mathcal T_{0,T}^0$ define the concatenation operator from $\Omega\times\Omega$ to $\Omega$ by the formula
$$ (\omega\otimes_\theta\omega')(u)=\omega(u)I_{[0,\theta(\omega)]}(u)+[\omega'(u)-\omega'(\theta(\omega))+\omega(\theta(\omega))] I_{(\theta(\omega),T]}(u).$$ 
Furthermore, for any function $\xi$ on $\Omega$ define the shifted function
$$ \xi^{\theta,\omega}(\omega')=\xi(\omega\otimes_\theta\omega'),\ \ \ \omega'\in\Omega.$$
We use the following properties of these objects (\cite{BayHua10}, Proposition A.1(ii) and Remark 2.5(ii)). 
\begin{lemma} \label{lem:5}
\begin{itemize}  
\item[(i)] Let $\xi$ be an integrable or nonnegative $\mathscr F_T$-measurable random variable and $\theta\in\mathcal T$. Then
$$ \mathsf E(\xi|\mathscr F_\theta)(\omega)=\int_\Omega \xi^{\theta,\omega}(\omega') \,d\mathsf P(\omega')\ \ \text{for}\ \mathsf P\text{-a.e.}\ \omega\in\Omega.$$
\item[(ii)] Let $\alpha\in\mathscr A$, $\theta\in\mathcal T_{t,T}^t$. Then for $\mathsf P$-a.e. $\omega\in\Omega$ we have
$$ X_u^{t,x,\alpha}(\omega\otimes_\theta\omega')=X_u^{\theta(\omega),X_\theta^{t,x,\alpha}(\omega),\alpha^{\theta,\omega}}(\omega'), \ \ u\in [\theta(\omega),T]\ \ \text{for}\ \mathsf P\text{-a.e.}\ \omega'\in\Omega.$$
\end{itemize}
\end{lemma}

Note, that although the assertion (ii) of this lemma is formulated in \cite{BayHua10} for a nonrandom $\theta=s$, it is used in the above form.  

\section{Proof of Theorem 1}
\setcounter{equation}{0}
\emph{Step 1.} Take $(t,x)\in[0,T]\times\mathbb R^d$, $\alpha\in\mathscr A$ and $\theta\in\mathcal T_{t,\tau^{t,x,\alpha}}^t$. We have
\begin{equation} \label{eq:3.1}
J(t,x,\alpha)=\mathsf E\left[\int_t^\theta f(s,X_s^{t,x,\alpha},\alpha_s)\,ds+\mathsf E\left(\int_\theta^{\tau^{t,x,\alpha}}f(s,X_s^{t,x,\alpha},\alpha_s)\,ds\biggr|\mathscr F_\theta\right)\right]
\end{equation}
Using Lemma \ref{lem:5}(ii), we get
$$\tau^{t,x,\alpha}(\omega\otimes_\theta\omega')=\tau_G(X_\cdot^{t,x,\alpha}(\omega\otimes_\theta\omega'))=
\tau_G\left(X_\cdot^{\theta(\omega),X_\theta^{t,x,\alpha}(\omega),\alpha^{\theta,\omega}}(\omega')\right)=\tau^{\theta(\omega),X_\theta^{t,x,\alpha}(\omega),\alpha^{\theta,\omega}}(\omega').$$
Furthermore, since $(\omega\otimes_\theta\omega')(t\wedge\theta(\omega))=\omega(t\wedge\theta(\omega))$ and $\theta(\omega(\cdot\wedge\theta(\omega))=\theta(\omega(\cdot))$ (see \cite{RevYor99}, Example 4.21.3$^\circ$), it follows that
$$\theta(\omega\otimes_\theta\omega')=\theta(\omega).$$
Hence, by Lemma \ref{lem:5}(i),
\begin{align} 
\mathsf E\left(\int_\theta^{\tau^{t,x,\alpha}}f(s,X_s^{t,x,\alpha},\alpha_s)\,ds\biggr|\mathscr F_\theta\right)(\omega)\nonumber\\
=\int_\Omega\int_{\theta(\omega)}^{\tau^{\theta(\omega),X_\theta^{t,x,\alpha}(\omega),\alpha^{\theta,\omega}}(\omega')} f(s,X_s^{\theta(\omega),X_\theta^{t,x,\alpha}(\omega),\alpha^{\theta,\omega}}(\omega'),\alpha^{\theta,\omega}(\omega'))\,ds\,d\mathsf P(\omega')\nonumber\\
=J(\theta(\omega),X_\theta^{t,x,\alpha}(\omega),\alpha^{\theta,\omega}) \label{eq:3.2}.
\end{align}

Put $\theta=t$. As $\alpha^{t,\omega}\in\mathscr A^t$, from (\ref{eq:3.1}), (\ref{eq:3.2}) we see that 
$$J(t,x,\alpha)=\int_\Omega J(t,x,\alpha^{t,\omega})\,d\mathsf P(\omega)\le v(t,x).$$
Thus, it is possible to pass from $\mathscr A^t$ to $\mathscr A$ in the definition of the value function:
\begin{equation} \label{eq:3.3}
 v(t,x)=\sup_{\alpha\in\mathscr A} J(t,x,\alpha). 
\end{equation} 
The above argumentation is completely analogous to \cite{BouTou11} (Remark 5.2). The representation (\ref{eq:3.3}) and Lemma \ref{lem:4} imply that $v$ is lower semicontinuous.

Now take $\alpha\in\mathscr A^t$ and  $\theta\in\mathcal T_{t,\tau^{t,x,\alpha}}^t$. As $\alpha^{\theta,\omega}\in\mathscr A^{\theta(\omega)}$, by formulas (\ref{eq:3.1}), (\ref{eq:3.2}) and the definition of $v$ it follows that
$$J(t,x,\alpha)\le\mathsf E\left[\int_t^\theta f(s,X_s^{t,x,\alpha},\alpha_s)\,ds+v(\theta,X_\theta^{t,x,\alpha})\right].$$
Evidently, this yields that
\begin{equation} \label{eq:3.4}
v(t,x) \le \sup_{\alpha\in\mathscr A^t}\inf_{\theta\in \mathcal T^t_{t,\tau^{t,x,\alpha}}}
\mathsf E\left[\int_t^\theta f(s,X_s^{t,x,\alpha},\alpha_s)\,ds+v(\theta,X_\theta^{t,x,\alpha})\right].
\end{equation}

\emph{Step 2. }Put $B(t,x;r)=\{(t',x')\in [0,T]\times\mathbb R^d: t'\in(t-r,t],\ |x'-x|<r\}$ and fix $\varepsilon>0$. For each $(t,x)\in[0,T]\times\mathbb R^d$ take an $\varepsilon$-optimal control $\alpha^{t,x}\in\mathscr A^t$:
$v(t,x)-\varepsilon\le J(t,x,\alpha^{t,x})$. For $t\in [0,T]$ and a nonnegative continuous function $\varphi\le v$ there exist  $r^{t,x}>0$ such that
$$ \varphi(t,x)+\varepsilon\ge \varphi(t',x'), \ \ \ J(t,x,\alpha^{t,x})-\varepsilon\le J(t',x',\alpha^{t,x})$$ 
for $(t',x')\in B(t,x;r^{t,x})$ since $(t',x')\mapsto J(t',x',\alpha^{t,x})$ is lower semicontinuous by Lemma \ref{lem:4}.
It follows that
\begin{equation} \label{eq:3.5}
J(t',x',\alpha^{t,x})\ge J(t,x,\alpha^{t,x})-\varepsilon\ge v(t,x)-2\varepsilon\ge\varphi(t,x)-2\varepsilon\ge\varphi(t',x')-3\varepsilon
\end{equation}
for $(t',x')\in B(t,x;r^{t,x})$.

Consider on $\mathbb R$ the upper limit topology with the basis $(a,b]$, $a,b\in\mathbb R$. Under this topology $\mathbb R$ is Lindel\"{o}f: see \cite{Mun00} (Chapter 4, \S30, Example 3) for the similar case of the lower limit topology. Since $\mathbb R^d$ under the standard topology is $\sigma$-compact, the space $\mathbb R\times\mathbb R^d$ is Lindel\"{of} under the product topology (see: \cite{AurTal12}, Fig.1) with the basis formed by the sets 
$$\{(t',x')\in \mathbb R\times\mathbb R^d: t'\in(t-r,t],\ |x'-x|<r\}.$$
In this topology the set $[0,T]\times\mathbb R^d$ is closed and hence, it is Lindel\"{o}f in the subspace topology with the basis $B(t,x;r)$, $(t,x)\in[0,T]\times\mathbb R^d$. 

Thus, we can take a countable cover $(B(t_i,x_i;r_i))_{i\in\mathbb N}$, $r_i=r^{t_i,x_i}$ of the set $[0,T]\times\mathbb R^d$. The family of disjoint sets $(A_i)_{i\in\mathbb N}$, defined by
$$ A_1=B(t_1,x_1;r_1),\ \ \ A_i=B(t_i,x_i;r_i)\backslash\left(\bigcup_{j=1}^{i-1} A_j\right),\ \ i\ge 2,$$
also covers $[0,T]\times\mathbb R^d$. For $\alpha\in\mathscr A^t$ and $\theta\in\mathcal T_{t,\tau^{t,x,\alpha}}^t$ put
$$\beta_s=\alpha_s I_{[0,\theta)}(s)+\sum_{i\in\mathbb N}\alpha^i_s I_{A_i}(\theta,X_\theta^{t,x,\alpha})I_{[\theta,T]}(s),\ \ \text{where } \alpha^i=\alpha^{t_i,x_i}.$$
It is easy to see that $\beta$ is $\mathbb F^t$-progressively measurable. Note that $\alpha^i\in\mathscr A^{t_i}\subset\mathscr A^{\theta(\omega)}$ if $(\theta,X_\theta^{t,x,\alpha})(\omega)\in A_i$. By (\ref{eq:3.1}), (\ref{eq:3.2}) and the definition of $\beta$, we have
\begin{align*}
J(t,x,\beta) &=\mathsf E\left[\int_t^\theta f(s,X_s^{t,x,\beta},\beta_s)\,ds+J(\theta,X_\theta^{t,x,\alpha},\beta^{\theta,\cdot})\right]\\
&=\mathsf E\left[\int_t^\theta f(s,X_s^{t,x,\alpha},\alpha_s)\,ds+\sum_{i\in\mathbb N}J(\theta,X_\theta^{t,x,\alpha},\alpha^i) I_{A_i}(\theta,X_\theta^{t,x,\alpha})\right]. 
\end{align*}

Furthermore, by (\ref{eq:3.5}),
\begin{align*}
\sum_{i\in\mathbb N}J(\theta,X_\theta^{t,x,\alpha},\alpha^i) I_{A_i}(\theta,X_\theta^{t,x,\alpha})
&\ge \sum_{i\in\mathbb N} (\varphi(\theta,X_\theta^{t,x,\alpha})-3\varepsilon) I_{A_i}(\theta,X_\theta^{t,x,\alpha})\\
&=\varphi(\theta,X_\theta^{t,x,\alpha})-3\varepsilon.
\end{align*}
It follows that
$$ v(t,x)\ge J(t,x,\beta)\ge \mathsf E\left[\int_t^\theta f(s,X_s^{t,x,\alpha},\alpha_s)\,ds + \varphi(\theta,X_\theta^{t,x,\alpha})\right]-3\varepsilon.$$
Since $v$ is lower semicontinuous, there exists a non-decreasing sequence $(\varphi_n)_{n=1}^\infty$ of nonnegative continuous functions convergent pointwise to $v$ on $[0,T]\times\mathbb R^d$ (see \cite{Bee93}, Theorem 1.3.7). By the monotone convergence theorem we get the inequality
$$ v(t,x)\ge \mathsf E\left[\int_t^\theta f(s,X_s^{t,x,\alpha},\alpha_s)\,ds + v(\theta,X_\theta^{t,x,\alpha})\right]-3\varepsilon.$$
By the arbitrariness of $\theta$, $\alpha$ and $\varepsilon$ it follows that
$$ v(t,x)\ge \sup_{\alpha\in\mathscr A^t}\sup_{\theta\in \mathcal T^t_{t,\tau^{t,x,\alpha}}}\mathsf E\left[\int_t^\theta f(s,X_s^{t,x,\alpha},\alpha_s)\,ds + v(\theta,X_\theta^{t,x,\alpha})\right].$$
Together with (\ref{eq:3.4}) this proves the assertion of Theorem 1.

\end{document}